\documentclass[a4paper,12pt]{article}
\usepackage{amssymb}
\usepackage{amsmath}
\usepackage{amsthm}
\usepackage{tikz}
\usepackage{verbatim}
\usetikzlibrary{arrows}

\newtheorem{theo}{Theorem}

\usepackage{amsfonts}
\usepackage{mathtools}
\usepackage{fullpage}
\usepackage{cancel}
\usepackage{color}

\newcommand{\pt}{\mbox{\rm proj}_\theta}
\renewcommand{\dh}{\dim_{\rm {H}}}
\newcommand{\dpk}{\dim_{\rm {P}}}
\newcommand{\proj}{ \mbox{\rm proj}}

\newcommand{\rn}{\mathbb{R}^n}

\title{Sixty Years of Fractal Projections}
\author{Kenneth Falconer, Jonathan Fraser and Xiong Jin}
\date{}

\begin{document}
\maketitle
\begin{abstract}
Sixty years ago, John Marstrand published a paper which, among other things, relates the Hausdorff dimension of a plane set to the dimensions of its orthogonal projections onto lines. For many years, the paper attracted very little attention. However, over the past 30 years, Marstrand's projection theorems have become the prototype for  many results in fractal geometry with numerous variants and applications and they continue to motivate leading research.
\end{abstract}

\section{Marstrand's 1954 paper}
\setcounter{equation}{0}
\setcounter{theo}{0}
\setcounter{figure}{0}
In 1954, John Marstrand's paper \cite{Mar} `Some fundamental geometrical properties of plane sets of fractional dimensions' was published in the Proceedings of the London Mathematical Society. The paper was essentially the work for his doctoral  thesis at Oxford, which was heavily influenced by Abram Besicovitch, a Russian born mathematician who pioneered geometric measure theory. For 25 years after its publication the paper attracted very limited attention, since then it has become one of the most frequently cited papers in the area now referred to as {\it fractal geometry}. Indeed, the paper was the first to consider the geometric properties of fractal dimensions.

The best-known results from the paper are the following two Projection Theorems, stated below in Marstrand's wording, which relate the dimensions of sets in the plane to those of their orthogonal projections onto lines through the origin. Note that `dimension' refers to Hausdorff dimension, and an `$s$-set' is a set that is measurable and of positive finite measure with respect to $s$-dimensional Hausdorff measure ${\mathcal H}^s$. `Almost all directions' means all lines making angle $\theta$ with the $x$-axis except for a set of $\theta \in [0,\pi)$ of Lebesgue measure 0.
\medskip

\noindent {\sc  Theorem I}. {\it Any $s$-set whose dimension is greater than unity projects into a set of positive Lebesgue measure in almost all directions.}
\medskip

\noindent {\sc  Theorem II}. {\it Any $s$-set whose dimension does not exceed unity projects into a set of dimension $s$ in almost all directions.}
\medskip

The statements are followed by a remark that, by a result of Roy Davies \cite{Dav}, every Borel or analytic set of infinite $s$-dimensional Hausdorff measure contains an $s$-set. This allows the theorems to be expressed in terms of Hausdorff dimension, and this is the form in which they are now usually stated. We write $\dh$ for Hausdorff dimension, ${\mathcal L}$ for Lebesgue measure, and $\pt$ for orthogonal projection of a set onto the line at angle $\theta$ to the $x$-axis, see Figure 1.

\begin{theo}{\rm \cite{Mar}} \label{marthm}
Let $E \subset \mathbb{R}^2$ be a Borel or analytic set. Then, for almost all $\theta \in [0,\pi)$,

(i) $\dh \mbox{\rm proj}_\theta E = \min\{\dh E, 1\}$,
\smallskip

(ii) ${\mathcal L} (\mbox{\rm proj}_\theta E) > 0$ if $\dh E >1$.
\end{theo}

 \begin{figure}[h]
\begin{center}
\includegraphics[scale=0.55]{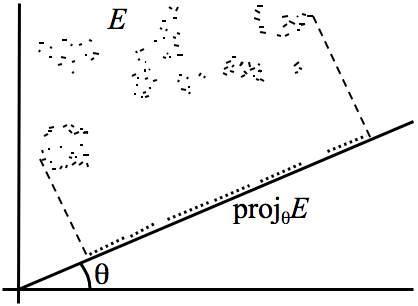}

Figure 1: Projection of a set $E$ onto a line in direction $\theta$.
\end{center}
\end{figure}

Since projection does not increase distances between points it follows easily from the definition of Hausdorff measure and dimension that $\dh \pt E \leq \min\{\dh E, 1\}$ for all $\theta$, but the opposite inequality is much more intricate. Marstrand's proofs depend heavily on plane geometry and measure theory, with, for example, careful estimates of the measures of narrow strips in various directions. As John Marstrand once remarked, analysis essentially consists of integrating functions in different ways and applying Fubini's theorem - but it may be difficult to find an appropriate function. The proofs in this paper illustrate this well. 

It is worth mentioning that Marstrand's paper \cite{Mar} includes a nice, but often forgotten, extension to the theorems,  that the same exceptional directions can apply to subsets of the given $s$-set that are of positive measure. In the following statement from the paper $|\quad|$ denotes Lebesgue measure.
\medskip

\noindent {\sc  Lemma 13}. {\it If E is an s-set and $s>1$, then for almost all angles $\theta$, all s-sets A which are contained in E satisfy}
$|{\rm proj}_\theta A|>0$.  
\medskip

Although Marstrand's paper is most often cited for the projection theorems, its 46 pages contain a great deal more, much of which anticipated other directions in fractal geometry.

$\bullet$ {\it Dimension of the intersection of sets with lines.} E.g.  Almost every line through ${\mathcal H}^s$-almost every point of an $s$-set $E$ ($s>1$) intersects $E$ in a set of dimension  $s-1$ and finite $s-1$-dimensional measure.
\smallskip

$\bullet$ {\it Construction of examples with particular projection properties.} E.g. For $1<s<2$ there exists an $s$-set which projects onto a set of  dimension $s-1$ in continuum many directions in every sector.
\smallskip

$\bullet$ {\it Dimension of exceptional sets.} The dimension of the set of points from which an irregular 1-set (see Section 4) has projection of positive Lebesgue measure is at most 1. 
\smallskip

$\bullet$ {\it Densities of $s$-sets.} The density $\lim_{r\to0} {\cal H}^s (E\cap B(x,r)) / (2r)^s$ of an $s$-set $E\subset \mathbb{R}^2$ can exist  and equal 1 for ${\cal H}^s$-almost all $x$ only if $s=0,1$ or 2. ($B(x,r)$ denotes the disc of centre $x$ and radius $r$.)
\smallskip

$\bullet$ {\it Angular densities.} Bounds are given for densities defined in segments emanating from points of an $s$-set.
\smallskip

$\bullet$ {\it Weak tangents.} For $1<s<2$ an $s$-set fails to have a  weak tangent (with an appropriate definition) almost everywhere.
\smallskip

This area of research is a central part of what is now termed fractal geometry.
This paper will survey the vast range of mathematics related to projections of sets that has developed over the past 60 years and which might be regarded as having its genesis in Marstrand's 1954 paper.

\section{The potential-theoretic approach}
\setcounter{equation}{0}
\setcounter{theo}{0}
\setcounter{figure}{0}

By virtue of the fact that an orthogonal projection is a Lipschitz map, we invariably have $\dh \proj E \leq \min\{m, \dh E\}$ for every set $E \subset \mathbb{R}^n$ and projection $\proj:\mathbb{R}^n\to V$ onto every $m$-dimensional subspace $V$, a fact that should be borne in mind throughout this article. It is inequalities in the opposite direction that require more work to establish. (However, a particularly straightforward situation is that for a connected set $E \subset \mathbb{R}^2$, both $\dh E \geq 1$ and  $\dh \proj E =1$ for projections onto lines in all directions with at most one exception.) Throughout this article we will always assume that the sets $E$ being projected are Borel or analytic -- pathological constructions show that dimension ceases to have useful geometric properties if completely general sets are considered.

Marstrand's proofs of his projection theorems were geometrically complicated and not particularly conducive to   extension or generalization.  
But in 1968  Kaufman \cite{Kau} gave new proofs of Theorem 1.1(i) using potential theory and of Theorem 1.1(ii)  using a Fourier transform method. This provided a rather more accessible approach, leading eventually to many generalizations and extensions. Kaufman's proofs depend on the following characterization of Hausdorff dimension in terms of an energy integral.
\begin{align}
\dh E = \sup\Big\{s : E & \mbox{ supports  a positive finite measure}\nonumber\\
& \mbox{ $\mu$ such that }
\int\int \frac{d\mu(x)d\mu(y)}{|x-y|^s} < \infty\Big\}.\label{energy}
\end{align}

Thus if $E \subset \mathbb{R}^2$ and $s<\dh E$ where $0< s<1$, we may find a measure $\mu$ supported by $E$ such that $ \displaystyle{\int \int\frac{d\mu(x)d\mu(y)}{|x-y|^s} < \infty}$.
Write $\mu_\theta$ for the projection of $\mu$ onto the line in direction $\theta$, so 
$\int_{-\infty}^\infty f(t)d\mu_\theta(t) = \int_E f(x\cdot  \theta) d\mu(x)$ for continuous $f$, where we identify $\theta$ with a unit vector in the direction $\theta$.  Then
\begin{eqnarray}
\int_{0}^\pi\bigg[ \int_{-\infty}^\infty \int_{-\infty}^\infty
\frac{d\mu_\theta(t)d\mu_\theta(u)}{|t-u|^s} \bigg]d\theta
&=&
\int_{0}^\pi\bigg[ \int_E \int_E
\frac{d\mu(x)d\mu(y)}{|x\cdot \theta-y\cdot \theta|^s} \bigg]d\theta\label{angint}\\
&=&
\int_E \int_E\int_{0}^\pi\frac{d\theta}{|{u}_{x-y}\cdot\theta|^s}  \frac{d\mu(x)d\mu(y)}{|x-y|^s}\nonumber\\
&\leq&c\int_E \int_E\frac{d\mu(x)d\mu(y)}{|x-y|^s}<\infty\nonumber
\end{eqnarray}
where $u_w$ denotes the unit vector $w/|w|$ and $\int_{0}^\pi |u_{x-y}\cdot\theta|^{-s}{d\theta}=c<\infty$. 

Hence for almost all $\theta$,  $ \displaystyle{\int_{-\infty}^\infty \int_{-\infty}^\infty
\frac{d\mu_\theta(t)d\mu_\theta(u)}{|t-u|^s} <\infty}$, so, since $\mu_\theta$ is supported by $\mbox{proj}_\theta E$, we conclude from the characterization \eqref{energy} that $\dh \mbox{proj}_\theta E\geq s$. This is true for all $s <\dh E$, so $\dh \mbox{proj}_\theta E\geq \dh E$ for almost all $\theta$.

For the case where $1<s<2$,  a variant of this argument  shows that
$$\int_{0}^\pi\bigg[ \int_{-\infty}^\infty \int_{-\infty}^\infty
|\widehat{\mu_\theta}(u)| ^2\bigg]du< \infty$$
where $\widehat{\mu_\theta}$ is the Fourier transform of $\mu_\theta$ from which it follows that $\mu_\theta$ is absolutely continuous with respect to Lebesgue measure with $L^2$ density, so in particular has support of positive Lebesgue measure.

In 1975 Mattila \cite{Mat} used potential theoretic methods to obtain the natural extension of these theorems to projections from higher dimensional spaces to subspaces. For $1\leq m <n$, and  $V$ an $m$-dimensional subspace of $\mathbb{R}^n$, let  $\mbox{proj}_V : \mathbb{R}^n\to V$. These subspaces form the Grassmanian $G(n,m)$, an $m(n-m)$-dimensional compact manifold which carries a natural invariant measure, locally equivalent to $m(n-m)$-dimensional Lebesgue measure. 
\begin{theo}{\rm \cite{Mat}}\label{higherdim}
Let $E \subset \mathbb{R}^n$ be a Borel or analytic set. Then, for almost all $V \in G(n,m)$,

(i) $\dh\mbox{\rm proj}_V E = \min\{\dh E, m\}$.
\smallskip

(ii) ${\mathcal L}^m (\mbox{\rm proj}_V E) > 0$ if $\dh E >m$, where ${\mathcal L}^m$ denotes $m$-dimensional Lebesgue measure on $V$.
\end{theo}

\section{Exceptional sets of projections}
\setcounter{equation}{0}
\setcounter{theo}{0}
\setcounter{figure}{0}

We can deduce rather more from Kaufman's proof above. Let $E \subset \mathbb{R}^2$ and $0<s<\dh E<1$. Let  $T=\{\theta : \dh \mbox{proj}_\theta E<s\}$. If $\dh T > s$ then it can be shown that we may find a measure $\nu$ supported by $T$ such that
$\int_T  |{\bf u}\cdot\theta|^{-s}d\nu(\theta)\leq c<\infty$ for every unit vector ${\bf u}$. If we integrate with respect to $\nu$ instead of Lebesgue measure in \eqref{angint} we still get a finite triple integral, and so for $\nu$-almost all $\theta\in T$ the $s$-energy of $\mu_\theta$ is finite and $\dh \mbox{proj}_\theta E\geq s$, a contradiction. It follows that 
if $E \subset \mathbb{R}^2$ and $0\leq s< \dh E < 1$ then
%$$\dh\{\theta : \dh \mbox{proj}_\theta E< \dh E\} \leq \dh E.$$
$$\dh\{\theta : \dh \mbox{proj}_\theta E< s\} \leq  s.$$

Thus the set of $\theta$ for which the projections have much smaller dimension than that of the set is correspondingly small. 
Indeed, the dimension of a projection is rarely less than half that of the set. As Bourgain \cite{Bou} and Oberlin \cite{Obe} showed, again when $E \subset \mathbb{R}^2$ and $ \dh E < 1$, 
$$\dh\{\theta : \dh \mbox{proj}_\theta E< \textstyle{\frac{1}{2}}\dh E\} =0.$$

For  $E \subset \mathbb{R}^2$ and $ \dh E > 1$,  the greater the `excess dimension' $ \dh E - 1$ the smaller the set of $\theta$ where Marstrand's conclusion fails. To be more precise:
$$\dh\{\theta : {\mathcal L}(\mbox{proj}_\theta E)=0\} \leq 2-\dh E.$$
This was first proved in \cite{Fal1} and all known proofs depend on Fourier transforms. 

Not surprisingly there are higher dimensional analogues of these bounds on the dimensions of the exceptional set, that is the set of $V \in G(n,m)$ for which the conclusions of Theorem \ref{higherdim} fail. 
These are summerised in the following inequalities, written for comparison with $m(n-m)$, the dimension of the Grassmanian $G(n,m)$, see \cite{Mat,Mat1,Mat3} for more details.

\begin{theo}\label{exceptproj}
Let $E\subset \mathbb{R}^n$ be a Borel or analytic set. 

\noindent (i) If $0<s < \dh E \leq m$ then
$$\dh\{V \in G(n,m) : \dh \mbox{\rm proj}_V E< s\} \leq m(n-m) -(m-s);$$
(ii) if $ \dh E \geq m$ then
$$\dh\{V \in G(n,m) : \dh \mbox{\rm proj}_V E< s\} \leq m(n-m) -(\dh E-s);$$
(iii) if $\dh E > m$  then
$$\dh\{V \in G(n,m) : {\mathcal L}^m(\mbox{\rm proj}_V E)=0\} \leq m(n-m)- (\dh E-m);$$
(iv) if $\dh E > 2m$  then
$$\dh\{V \in G(n,m) : \mbox{\rm proj}_V E \mbox{ \rm has  empty interior} \}\leq m(n-m)-  (\dh E-2m).$$
\end{theo}

\section{Sets of integer dimension}
\setcounter{equation}{0}
\setcounter{theo}{0}
\setcounter{figure}{0}

Marstrand was the first person to consider the effect of projection on the numerical value of the dimension, but his paper also includes a few results on
projections of $s$-sets in the `critical case' where $s=1$.  This case had been studied in great detail somewhat earlier by Besicovitch around the 1930s \cite{Bes1,Bes2,Bes3} who showed that 1-sets or `linearly-measurable sets' in the plane could be decomposed into a regular part and an irregular part, defined in terms of local densities  $D(x) = \lim_{r\to0} {\cal H}^1 (E\cap B(x,r)) / 2r$. The {\it regular}  part consists of those $x$ where the limit $D(x)$ exists with $D(x)=1$, and the  {\it irregular}  part is formed by the remaining points. Besicovitch showed that, to within a set of measure 0,  the regular part is `curve-like', that is a subset of a countable collection of rectifiable curves. On the other hand, the irregular part is `dust-like' intersecting every rectifiable curve in length 0. Using intricate geometrical arguments, Besicovitch obtained the following projection theorem.

\begin{theo}{\rm \cite{Bes3}} \label{besproj}
Let $E \subset \mathbb{R}^2$ be a $1$-set. 

(i) If $E$ is regular then ${\mathcal L} (\mbox{\rm proj}_\theta E) > 0$ for all $\theta \in [0,\pi)$ except perhaps for a single value of $\theta$.

(ii)  If $E$ is irregular then ${\mathcal L} (\mbox{\rm proj}_\theta E) = 0$ for almost all $\theta \in [0,\pi)$.
\end{theo}

The natural higher dimensional versions of Theorem \ref{besproj}, with appropriate definitions of regular and irregular sets, were obtained by Federer \cite{Fed, Fed1}.

If $E$ is measurable and of $\sigma$-finite $ {\cal H}^1$ measure, it follows from Theorem \ref{besproj}  that  ${\mathcal L} (\mbox{\rm proj}_\theta E)$ is either 0 for almost all $\theta$ or positive  for almost all $\theta$, by decomposing $E$  into countably many 1-sets. However if  $\dh E =1$ but $E$ is not $\sigma$-finite then strange things can occur: we can find a set $E$ whose projections  are, to within Lebesgue measure 0, anything we like.

\begin{theo}{\rm \cite{Dav,Fal2}} \label{digsun}
For each $\theta \in [0,\pi)$ let $E_\theta$ be a given subset of the line through the origin of $\mathbb{R}^2$ in direction $\theta$, such that $\bigcup_{0\leq \theta <\pi} E_\theta$ is plane Lebesgue measurable. Then there exists a Borel set $E \subset \mathbb{R}^2$ such that,  for almost all $\theta$, ${\mathcal L} (E_\theta \,\triangle\, \proj_\theta E)  =0$ where $\triangle$ denotes symmetric difference, in other words $\proj_\theta E$ differs from the prescribed set $E_\theta$ by a set of negligible length.
\end{theo}
Theorem \ref{digsun} may be obtained by dualising a result of Davies \cite{Dav} on covering a plane set by lines without increasing its plane Lebesgue measure. It was proved directly, along with the natural higher dimension analogues, in \cite{Fal2}. For projections from $\mathbb{R}^3$ to $\mathbb{R}^2$ this has become known as the `digital sundial theorem':  Given a subset  $E_V$  of  each 2-dimensional subspace $V$ of $\mathbb{R}^3$ (with a measurability condition), there exists a Borel set $E \subset \mathbb{R}^3$ such that, for almost all subspaces $V$,
${\cal L}^2(E_V \triangle{\rm proj}_V E) = 0$. Thus, in theory at least, there is a set in space such that the shadow cast by the sun gives the thickened digits of the time at any instant.

\section{Packing dimensions}\label{pacsec}
\setcounter{equation}{0}
\setcounter{theo}{0}
\setcounter{figure}{0}

Packing measures and packing dimension were  introduced  by Tricot \cite{Tri}  in 1982 as a sort of dual to their Hausdorff counterparts, see \cite{Fal,Mat1}. Whilst packing measures require an extra step in their definition, the gap of over sixty years between the two concepts seems very surprising with hindsight. Nowadays, however, every problem that  involves Hausdorff dimension is almost routinely studied in terms of packing dimension as well. Projection theorems are no exception, but the dimensional relationships turn out to be more complicated in the packing dimension case. 

J\"{a}rvenp\"{a}\"{a} \cite{Jar} constructed compact sets $E\subset \mathbb{R}^n$ with $\dpk E$ taking any prescribed value in $(0,n]$ such that
$ \dpk \mbox{proj}_V E =  \dpk E\big/\big(1+ (1/m - 1/n)  \dpk E\big)$  for all $V\in G(n,m)$. 
This is essentially the least value that can be obtained, that is 
$$ \frac{\dpk E} {1+ (1/m - 1/n)  \dpk E} \leq  \dpk \mbox{proj}_V E \leq \min\{\dpk E,m\}$$
for almost all $V\in G(n,m)$, see \cite{FH}. For packing dimensions of projections of measures, rather than sets, this lower bound was refined to incorporate both the Hausdorff and packing dimensions of the measure, see \cite{FM}.
  
These inequalities  raised the question of whether $\dpk \mbox{proj}_V E$ takes a common value for almost all subspaces $V$ and this was answered affirmatively with the introduction of `dimension profiles' \cite{FH2}. The   packing dimension profile $\dpk^s E$ of a set $E\subset \mathbb{R}^n$ reflects how $E$ appears when viewed in an $s$-dimensional setting. 
For $s>0$ the $s$-{\it dimensional packing dimension profile} of a measure $\mu$ on $\mathbb{R}^n$ with bounded support is defined in terms of local densities of measures with respect to a kernel of the form $\min\{1,r^s/|x-y|^s\}$:
$$\dpk^s\mu  = \sup \bigg\{t \geq 0: \liminf_{r \searrow 0} r^{-t} \int \min\Big\{1,\frac{r^s}{|x-y|^s}\Big\}d\mu(y)<\infty 
\mbox{ for } \mu \mbox{-almost all } x \in \rn\bigg\}.$$
This leads to the $s$-{\it dimensional packing dimension profile} of a set $E\subset \mathbb{R} ^n$
$$\dpk^s E  = \sup \big\{\dpk^s \mu: \mu \mbox{ is a finite compactly supported measure on } E\},$$ 
see \cite{FH2}. The profiles generalize packing dimensions, since $\dpk^n E = \dpk E$ for $E \subset \rn$. The profiles may also be expressed in terms of measures defined by weighted coverings, see \cite{How,KX}.
\begin{theo}{\rm \cite{FH2}}\label{packdim}
Let $E \subset \mathbb{R}^n$ be a Borel or analytic set. Then, for almost all $V \in G(n,m)$,
$$\dpk \mbox{\rm proj}_V E = \dpk^m E.$$
\end{theo}
There is a certain parallel with Hausdorff dimensions, where one might define a dimension profile simply as 
$\dh^s E = \min\{s,\dh E\}$ which, by Marstrand's theorem, gives the almost sure Hausdorff dimension of projections onto $s$-dimensional subspaces. 

As well as giving the almost sure packing dimension of the projections, the profiles provide upper bounds for the dimension of the exceptional set of directions for which the packing dimension falls below the almost sure value.

Since their introduction, packing dimension profiles have cropped up in other contexts, notably to give the almost sure packing dimension of images of sets under fractional Brownian motion \cite{KX,Xia}.

\section{Projections in restricted directions}
\setcounter{equation}{0}
\setcounter{theo}{0}
\setcounter{figure}{0}
A general question that has been around for many years is under what circumstances we can get projection theorems for projections onto families of lines or subspaces that form proper subsets of $V(n,m)$.  For instance, if
$\{\theta(t): t \in P\}$ is a smooth curve or submanifold 
 of $V(n,m)$ smoothly  parameterized by a set $P \subset \mathbb{R}^k$, then what can we conclude about
 $\dh \proj_{\theta(t)} E$ for ${\mathcal L}^k$-almost all $t \in P$, where  ${\mathcal L}^k$ is $k$-dimensional Lebsegue measure?
 
 For a simple example, it follows easily from Theorem \ref{exceptproj} (ii)-(iii) that if $\{\theta(t): 0\leq t \leq 1\}$ is a smoothly parameterized curve of directions in $\mathbb{R}^3$ (i.e. a curve in $V(3,1))$, then for almost all $0 \leq t \leq 1$ we have $\dh \proj_\theta E\geq \min\{\dh E -1,1\}$ and if $\dh E > 2$ then
 $ {\mathcal L}^1(\proj_{\theta (t)} E)>0$, where $\proj_{\theta (t)}$ denotes projection onto the line in direction $\theta(t)$.

The following lower bounds  were obtained J\"{a}rvenp\"{a}\"{a}, J\"{a}rvenp\"{a}\"{a} and Keleti \cite{JJK} for parameterized families of projections from  $\mathbb{R}^n$ to $m$-dimensional subspaces, see also \cite{JJLL}. 
  For $0<k<m(n-m)$ define the integers
$$p(l) = n-m -\bigg\lfloor \frac{k-l(n-m)}{m-l}\bigg\rfloor \qquad (l = 0,1,\ldots, m-1),$$
where the `floor' symbol `$\lfloor x\rfloor$' denotes the largest integer no greater than $x$.

\begin{theo}{\rm \cite{JJK}}\label{projgeneral}
Let $P \subset \mathbb{R}^k$ be an open parameter set and let $E \subset \mathbb{R}^n$ be a Borel or analytic set. Let $\{V(t)\subset G(n,m): t \in P\}$ be a  family of subspaces such that $V$ is $C^1$ with the derivative $D_t V(t)$ injective for all $t \in P$. Then, for all $l=0,1,\ldots,m$ and ${\cal L}^k$-almost all $t \in P$,
$$
\dh\proj_{V(t)} E \geq  
\left\{
\begin{array}{lrl}
 \dh E-p(l) & \mbox{ if } \  p(l) + l  \!\! & \leq \dh E \leq p(l) + l +1  \\
l+1 & \mbox{ if }  \ p(l) + l +1 \!\! &\leq \dh E \leq p(l+1) + l +1
\end{array}
\right. .
$$
Moreover, if $\dh E > p(m-1) + m$ then  ${\cal L}^m ( \proj_{V(t)} E) > 0$ for ${\cal L}^k$-almost all $t \in P$.
\end{theo}

These are the best possible bounds for general parameterized families of projections.
The same paper \cite{JJK} includes generalizations of these results to smoothly parameterized families of $C^2$-mappings.

Better lower bounds may be obtained  if there is curvature in the mapping $s \mapsto V(s)$. This is a difficult area, and work to date mainly concerns projections from $\mathbb{R}^3$ to lines and planes. 
Let  $\theta: (0,1) \to S^2$ be a family of directions given by  a $C^3$-function $\theta$, where $S^2$ is the 2-sphere embedded in $\mathbb{R}^3$. We say that the curve of directions  is {\it non-degenerate} if 
$$\mbox{span }\{\theta(t), \theta'(t), \theta''(t)\} = \mathbb{R}^3 \quad \mbox{ for all  } t \in (0,1).$$

The following theorem was  proved by recently  by F\"{a}ssler and Orponen \cite{FO1}.

\begin{theo}{\rm \cite{FO1}}\label{projline}
Let $E \subset \mathbb{R}^3$ be a Borel or analytic set, let $\theta(t)$ be a non-degenerate family of directions, and let $\proj_{\theta(t)}$ denote projection onto the line in direction $\theta$. Then, for almost all $t \in (0,1)$,
\begin{equation}\label{projlineeq}
\dh\proj_{\theta(t)} E \geq \min\{\dh E, \textstyle{\frac{1}{2}}\}.
\end{equation}
\end{theo}
It is conjectured that $\frac{1}{2}$ can be replaced by $1$ in \eqref{projlineeq} and this is verified where $E$ is a self-similar sets without rotations in \cite{FO1}, a paper that also contains estimates for packing dimensions of projections.

The following bounds have been established for projections onto planes in $\mathbb{R}^3$ in the non-degenerate case. The conjectured lower bound is $\min\{ \dh E, 2\} $ and the bound  $\min\{ \dh E, 1\}$ for all values of $\dh E$ was obtained in \cite{FO1}. The further improvements stated below come from Fourier restriction methods \cite{OO}.
\begin{theo}{\rm \cite{FO1,OO}}\label{projplane}
Let $E \subset \mathbb{R}^3$ be a Borel or analytic set, let $\theta(t)$ be a non-degenerate family of directions, and let $\proj_{V_\theta(t)}$ denote projection onto the plane perpendicular to direction $\theta$. Then, for almost all $t \in (0,1)$,
\begin{equation}\label{projplaneeq}
\dh\proj_{V_\theta(t)} E \geq  
\left\{
\begin{array}{ll}
\min\{ \dh E, 1\} & \mbox{ if }   0 \leq \dh E \leq \frac{4}{3}  \\
 \textstyle{\frac{3}{4}} \dh E & \mbox{ if } \frac{4}{3} \leq \dh E \leq 2  \\
  \min\{ \dh E- \frac{1}{2}, 2\} &  \mbox{ if } 2 \leq \dh E \leq 3   
\end{array}
\right. .
\end{equation}
\end{theo}
Orponen \cite{Orp1} also showed that there exist numbers  $\sigma(\lambda)>1$ defined for $\lambda >1$, and increasing with $\lambda$, such that if $\dh E>1$ then $\dh\proj_{V_\theta(t)}E\geq \sigma(\dh E)$ for almost all $t$ .

Estimates for packing dimensions of projections may be found in \cite{FO1}.
The introduction of the paper \cite{Orp1} provides a recent overview of this area.

\section{Generalized  projections}
\setcounter{equation}{0} 
\setcounter{theo}{0}
\setcounter{figure}{0}

The projection theorems are a special case of much more general results. The essential property in Kaufman's proof in Section 2 is that the integral over the parameter $\theta$ satisfies
$\int |\proj_\theta x -\proj_\theta y|^{-s} d\theta \leq c |x-y|^{-s}$; such a condition can hold for many other parameterized families of mappings as well as for $\proj_\theta$.

Thus for $X \subset \mathbb{R}^n$ a compact domain, consider a family of maps
$\pi_\theta:  X\to  \mathbb{R}^m$ for $\theta$ in an open  parameter set $P \subset \mathbb{R}^k$.  Assume that the derivatives with respect to $\theta$, $D_\theta\pi_\theta (x)$  exist and are bounded. 

Let 
$$\Phi_\theta (x,y) = \frac{|\pi_\theta(x) - \pi_\theta(y)|}{|x-y|}.$$
The family $\{\pi_\theta: \theta\in P\}$ is {\it transversal} if there is a constant $c$ such that
\begin{equation}\label{transverse}
|\Phi_\theta (x,y)| \leq c \implies \det \big(D_\theta \Phi_\theta (x,y)(D_\theta \Phi_\theta (x,y))^T\big) \geq c
\end{equation}
for $\theta \in P$ and $x,y \in X$,
where $D_\theta$ denotes the derivative with respect to $\theta$ and $T$ denotes the transpose of a matrix.
This condition implies that if $\theta \in P$ is such that $\Phi_\theta (x,y)$ is small, then $\Phi_\theta (x,y)$ must be varying reasonably fast as $\theta$ changes in a direction perpendicular to the kernel of the derivative matrix.

By generalizing beyond recognition earlier arguments involving potential theory and Fourier transforms,  
Peres and Schlag \cite{PS} obtained theorems such as the following for a transversal family of generalized projections; compare Theorem \ref{exceptproj}.

\begin{theo}{\rm \cite{PS}}\label{genproj}
For $X \subset \mathbb{R}^n$  and $P \subset \mathbb{R}^k$, let $\{\pi_\theta:  X \to  \mathbb{R}^m: \theta \in P\}$ be a transversal family   and let $E\subset X$ be a Borel set. 

\noindent (i) If $0< t < \dh E \leq m$ then
$$\dh\{\theta \in P : \dh \pi_\theta E< t\} \leq k -(m-t),$$
(ii) if $ \dh E > m$ then
$$\dh\{\theta \in P : \dh \pi_\theta E< t\} \leq k -(\dh E-t),$$
(iii) if $\dh E > m$  then
$$\dh\{\theta \in P : {\mathcal L}^m(\pi_\theta E )=0\} \leq k -( \dh E-m),$$
(iv) if $\dh E > 2m$  then
$$\dh\{\theta \in P : \pi_\theta E \mbox{ \rm has  empty interior} \}\leq n-  \dh E+2.$$
\end{theo}

This powerful result has been applied to many situations,  including Bernoulli convolutions, sums of Cantor sets and pinned distance sets, see \cite{PS}. For a recent treatment of transversality, see \cite{Mat2}.

Leikas \cite{Lei} has used transversality to extend the packing dimension conclusions of Section \ref{pacsec}  to families of mappings between Riemannian manifolds where the dimension profiles again play a central role.

\section{Projections of self-similar and self-affine sets}
\setcounter{equation}{0}
\setcounter{theo}{0}
\setcounter{figure}{0}

One of the difficulties with the projection theorems is that they tell us nothing about the dimension or measure of the projection in any given direction.
There has been considerable recent interest in examining the dimensions of projections in specific directions for particular sets or classes of sets, and especially in finding sets for which the conclusions of Marstrand's theorems are valid for all, or virtually all, directions. Of particular interest are self-similar sets. 

Recall that an {\it iterated function system} (IFS)  is a family of  contractions $\{f_1, \ldots,f_k\}$ with $f_i : \mathbb{R}^n\to \mathbb{R}^n$. An IFS determines a unique non-empty compact  $E\subset \mathbb{R}^d$ such that 
\begin{equation}\label{attractor}
E= \bigcup_{i=1}^k f_i(E),
\end{equation}
called the {\it attractor} of the IFS, see \cite{Fal,Hut}. If the $f_i$ are all similarities, that is of the form
\begin{equation}\label{simform}
f_i(x) = r_i O_i (x) + a_i,
\end{equation}
where $0<r_i<1$ is the contraction ratio, $O_i$ is an orthonormal map, i.e. a rotation or reflection, and $a_i$ is a translation,
then $E$ is termed {\it self-similar}.
An IFS of similarities satisfies the {\it strong separation  condition} (SSC) if the union (\ref{attractor}) is disjoint, and the {\it open set condition} (OSC) if there is a non-empty open set $U$ such that $\cup_{i=1}^k f_i(U)\subset U$ with this union disjoint. If either SSC or OSC hold then $\dh E = s$, where $s$ is the {\it similarity dimension} given by $\sum_{i=1}^k r_i^s=1$, where $r_i$ is the similarity ratio of $f_i$, and moreover $0<{\mathcal H}^s(E) < \infty$. The {\it rotation group} $G= \langle O_1, \ldots, O_k \rangle$ generated by the orthonormal components of the similarities plays a crucial role in the behaviour of the projections of self-similar sets.

It is easy to construct self-similar sets with a finite rotation group $G$ for which the conclusions of Marstrand's theorem fail in certain directions. For example, let $f_1,\ldots,f_4$ be  homotheties (that is similarities with $O_i$ the identity in \eqref{simform}) of ratio $0<r<\frac{1}{4}$  that map the unit square $S$ into itself, each $f_i$ fixing one of the four corners. Then $\dh E = -\log 4 /\log r$,   but the projections of $E$  onto the sides of the square have  dimension $ -\log 2 /\log r$
 and onto the diagonals of $S$ have dimension $ -\log3 /\log r$, a consequence of the alignment of  the component squares $f_i(S)$ under projection. There is a similar reduction in the dimension of projections in directions $\theta$ whenever $ {\rm proj}_\theta (f_{i_1} \circ \cdots \circ f_{i_k}(S)) =  {\rm proj}_\theta (f_{j_1} \circ \cdots \circ f_{j_k}(S))$ for distinct words $i_1, \ldots, i_k$ and $j_1, \ldots, j_k$.
         
Kenyon \cite{Ken} conducted a detailed investigation of the projections of the 1-dimensional Sierpi\'{n}ski gasket $E\subset \mathbb{R}^2$, that is the self-similar set defined by the similarities
$$\textstyle{f_1(x,y) = (\frac{1}{3}x,\frac{1}{3}y),\  f_2(x,y) = (\frac{1}{3}x +\frac{2}{3},\frac{1}{3}y),
\ f_3(x,y) = (\frac{1}{3}x,\frac{1}{3}y+\frac{2}{3})}.$$
He showed that the projection of $E$ onto a line making an angle to the $x$-axis with tangent $p/q$ with  has dimension strictly less than 1 if $p+q \not\equiv 0 \, (\mbox{mod } 3)$, but if $p+q \equiv 0 \, (\mbox{mod } 3)$ then the projection has non-empty interior. For irrational directions he proved that the projections have Lebesgue measure 0 and  Hochman \cite{Hoc1} complemented this by showing that they nevertheless have Hausdorff dimension 1.

In fact, when the rotation group is finite, there are always some projections for which direct overlapping of the projection of components of the usual iterated construction leads to measure 0, as the following theorem of Farkas shows.
\begin{theo}{\rm \cite{Far}}\label{mesproj}
If $E\subset \mathbb{R}^n$ is self-similar with finite rotation group $G$ and similarity dimension $s$,  then 
$\dh\mbox{\rm proj}_V E<s $ for some $V \in G(n,m)$. In particular if $E$ satisfies OSC and $0< \dh E < m$ then 
$\dh\mbox{\rm proj}_V E<\dh E $ for some $V$.
\end{theo}

A rather different situation occurs if the IFS has {\it dense rotations}, that is the rotation group $G$ is dense in the full group of rotations $SO(n,\mathbb{R})$ or in the group of isometries 
$O(n,\mathbb{R})$. Note that an IFS of similarities of the plane has dense rotations if at least one of the rotations in the group is an irrational multiple of $\pi$. 

\begin{theo}{\rm \cite{PS,HS}}\label{denserots}
If $E\subset \mathbb{R}^n$ is self-similar with dense rotations   then 
\begin{equation}\label{hocshm}
\dh\mbox{\rm proj}_V E = \min\{\dh E, m\} \mbox{ for  all } V \in G(n,m).
\end{equation}
More generally, $\dh g(E) = \min\{\dh E, m\}$ for all $C^1$ mappings $g: E \to \mathbb{R}^m$ without singular  points, that is maps with non-singular derivative matrix.
\end{theo}
Peres and Shmerkin \cite{PS} proved \eqref{hocshm} in the plane without requiring any separation condition on the IFS. To show this they set up a discrete version of Marstrand's projection theorem to construct a tree of intervals in the subspace (line) $V$  followed by an application of Weyl's equidistribution theorem.  Hochman and Shmerkin \cite{HS}  proved the theorem in higher dimensions, including the extension to $C^1$ mappings, for $E$  satisfying the open set condition. Their proof uses the CP-chains of Furstenberg \cite{Fur,Fur1}, see also \cite{Hoc},  and has three main ingredients: the lower semicontinuity of the expected Hausdorff dimension of the projection of a measure with respect to its `micromeasures', Marstrand's projection theorem, and the invariance of the dimension of projections under the action of the rotation group.

That the open set condition is not essential follows since, for all $ \epsilon >0$,  we can use a Vitali argument to set up a new IFS, consisting of compositions of the $f_i$, that satisfies SSC, with attractor $E' \subset E$ such that $\dh E' > \dh E - \epsilon$; we can also ensure that the new IFS has dense rotations if the original one has, see \cite{FJ,Far,PS,Orp3}.

It is also natural to ask about  the Lebesgue measures of the projections of self-similar sets. We have seen examples of self-similar sets $E$ of Hausdorff dimension $s<m$ with finite rotation group and satisfying SSC such that $ {\mathcal H}^s (\mbox{\rm proj}_V E)$ is  positive for some subspaces $V$ and 0 for others.  For dense rotations, the situation is clear cut: the following theorem was proved by Ero\u{g}lu \cite{Ero} in the plane case for projections when OSC is satisfied, and for more general mappings with the separation condition removed by Farkas \cite{Far}.

\begin{theo}{\rm \cite{Far}}\label{projfar}
Let $E \subset \rn$ be the self-similar attractor of an IFS with dense rotations, with $\dh E = s$. Then ${\mathcal H}^s (\mbox{\rm proj}_V E) = 0$ for all $V \in G(n,m)$. More generally, ${\mathcal H}^s (g( E)) = 0$ for all $C^1$ mappings $g: E \to \rn$ without singular points.
\end{theo}

From Theorem \ref{denserots}, if $\dh E>m$ then in the dense rotation case $\dh \mbox{\rm proj}_V E = m$ for all $V \in G(n,m)$, but we might hope from the second part of Marstrand's theorem that the projections also have positive Lebesgue measure. Shmerkin and Solomyak showed that this is very nearly so in the plane.
 
 \begin{theo}{\rm \cite{SS}}\label{projSS}
Let $E \subset \mathbb{R}^2$ be the self-similar attractor of an IFS with dense rotations with $1<\dh E<2$. Then ${\mathcal L}^1 (\mbox{\rm proj}_\theta E) >0$ for all $\theta$ except for a set of $\theta$ of Hausdorff dimension 0. \end{theo}

The proof depends on careful estimation of the decay of the Fourier transforms of projections of a measure supported by $E$. The method can be traced back to  
a study of Bernoulli convolutions by Erd\H{o}s \cite{Erd}, which Kahane \cite{Kah} pointed out gave an exceptional set of Hausdorff dimension 0 rather than just Lebesgue measure 0, see \cite{PSS}.

The attractor of an IFS  is {\it self-affine} if \eqref{attractor} holds for affine contractions $\{f_1,\ldots,f_k\}$. A plane self-affine set is a  {\it carpet} if the contractions are of the form 
\begin{equation}
f_i(x,y) = (a_i x +c_i, b_i y +d_i), \label{carpetIFS}
\end{equation}
i.e. affine transformations that leave the horizontal and vertical directions invariant. For many  self-affine carpets the dimensions of the projections behave well except in directions parallel to the axes.

 \begin{figure}[h]
\begin{center}
\includegraphics[scale=0.34]{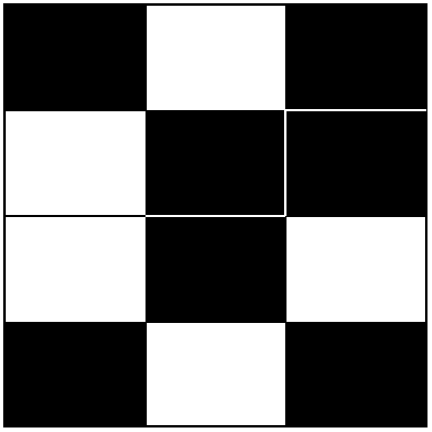}
\qquad
\includegraphics[scale=0.34]{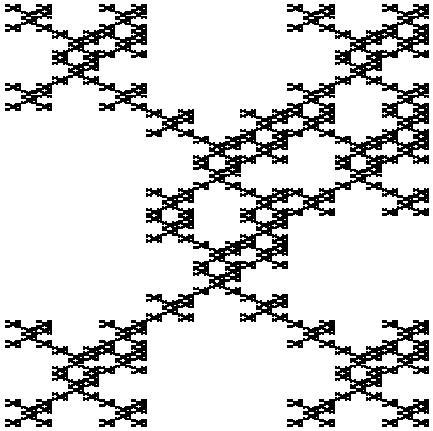}

Figure 2: A Bedford McMullen self-affine carpet obtained by repeated substitution of the left-hand pattern in itself.
\end{center}
\end{figure}

\begin{theo}{\rm \cite{FJS}}\label{projFJS}
Let $E \subset \mathbb{R}^2$ be a self-affine carpet in the Bedford-McMullen, Gatzouras-Lalley or Bara\'{n}ski class. If the IFS is of irrational type, then  
$\dh \proj_\theta E = \min\{\dh E,1\}$ for all $\theta$ except possibly $\theta= 0$ and  $\theta= \frac{1}{2}\pi.$
\end{theo}
For definitions and details of these different classes of carpets, see \cite{FJS}. The IFS is of {\it irrational type} if, roughly speaking, $\log a_i /\log b_i$ is irrational for at least one of the $f_i$ in \eqref{carpetIFS}.

Along similar lines, for an integer $n \geq 2$, let $T_{n}:[0,1] \to[0,1]$ (where 0 and 1 are identified) be given by $T_{n}(x) =nx\,(\mbox{\rm mod} 1)$.  In the 1960s Furstenberg conjectured that if $E$ and $F$ are closed sets invariant under $T_2$ and $T_3$ respectively, then $\dh \proj_\theta (E \times F)$ should equal 
$\min\{\dh (E \times F),1\}$ for all $\theta$ except possibly $\theta= 0$ and  $\theta= \frac{1}{2}\pi.$  This was proved by Hochman and Shmerkin \cite{HS} along with more general results such as the following.

\begin{theo}{\rm \cite{HS}}
Let $E$ and $F$ be closed subsets of $[0,1]$ invariant under $T_m$ and $T_n$ respectively, where $m,n$ are not powers of the same integer.  Then $\dh \proj_\theta (E \times F)=\min\{\dh (E \times F),1\}$ for all $\theta$ except possibly $\theta= 0$ and  $\theta= \frac{1}{2}\pi$.
\end{theo}

Projection properties of self-affine measures underpin this work and there are measure analogues of these theorems, see \cite{FFS,FJS,HS}.

\begin{comment}
also obtained the analogous result for measures.  We write $\dh \mu =\inf\{\dh A : \mu(A)>0\}$ to denote the {\it lower Hausdorff dimension of a measure} $\mu$.
Continuing this line of thought, for integers $m,n \geq 2$, let $T_{m,n}:[0,1] \times[0,1]\to[0,1] \times[0,1]$ (where 0 and 1 are identified) be given by $T_{m,n}(x,y) =(T_m(x), T_n(y))$. Products of sets or measures invariant under $T_m$ and $T_n$ respectively are easily seen to be invariant under $T_{m,n}$, but general $T_{m,n}$ invariant sets or measures can be much more complicated. Self-affine carpets in the Bedford-McMullen class are $T_{m,n}$-invariant, with the associated self-affine Bernoulli measures  $T_{m,n}$ invariant measures.  The following measure theoretic analogue of Theorem \ref{projFJS} was proved in \cite{FFS}, but the extension to general $T_{m,n}$ invariant sets and measures remains open.
 \begin{theo}{\rm \cite{FFS}}
Let $\mu$ be a self-affine measure for a carpet in the Bedford-McMullen class. Then $\dh \proj_\theta \mu = \min\{\dh \mu,1\}$ for all $\theta$ except possibly $\theta= 0$ and  $\theta= \frac{1}{2}\pi.$
\end{theo}
\end{comment}

\section{Projections of random sets}
\setcounter{equation}{0}
\setcounter{theo}{0}
\setcounter{figure}{0}

Fractal percolation provides a natural method of generating statistically self-similar fractals, with the same random process determining the form of the fractals at both small and large scales. 

Best known is Mandelbrot percolation, based on repeated decomposition of squares into  smaller subsquares from which a subset is selected at random.
Let $D$ denote the unit square in $\mathbb{R}^2$. 
Fix an integer $M \geq 2$ and  a probability $0<p<1$. We divide $D$ into $M^2$ closed subsquares of side $1/M$ in the obvious way, and retain each subsquare independently with probability $p$ to get a set $D_1$ formed as a union of the retained subsquares. We repeat this process with the squares of $D_1$, dividing each into $M^2$ subsquares of side $1/M^2$ and choosing each with probability $p$ to get a set $D_2$, and so on. This  leads to the random {\it percolation set} $E = \cap_{k=0}^\infty D_k$.
If $p > M^{-2}$ then there is a positive probability of non-extinction, i.e. that $E\neq \emptyset$, conditional on which
$\dh E = 2+\log p / \log M$ almost surely. 

The topological properties of Mandelbrot percolation have been studied  extensively, see \cite{Dek, Fal,RS3} for surveys. In particular there is a critical probability $p_c$ with $1/M < p_c<1$ such that if $p>p_c$ then, conditional on non-extinction, $E$ contains many connected components, so projections onto all lines  automatically have positive Lebesgue measure.  If $p \leq p_c$ the percolation set $E$ is totally disconnected, and Marstrand's theorems provide information on projections of $E$ in almost all directions. However, Rams and Simon \cite{RS,RS2, RS3} recently showed using a careful geometrical analysis that, conditional on $E\neq \emptyset$, almost surely the conclusions of Theorem \ref{marthm} hold for all projections. 

\begin{theo}{\rm \cite{RS}} \label{manperc}
Let $E$ be the random set obtained by the Mandelbrot percolation process in the plane based on subdivision of squares into $M^2$ subsquares, each square being retained with probability $p>1/M^2$. Then, with positive probability $E\neq \emptyset$, conditional on which:

(i) $\dh \mbox{\rm proj}_\theta E = \min\{\dh E, 1\}$ for all $\theta \in [0,\pi),$
\smallskip

(ii) if $p >1/M$ then for all $\theta \in [0,\pi)$, $\proj_\theta E$ contains an interval and in particular ${\mathcal L} (\mbox{\rm proj}_\theta E) > 0$.
\end{theo}

The natural higher dimensional analogues of this theorem for projections onto all $V \in G(n,m)$ are also valid, see \cite{SV}.
There are also versions of this result when the squares are selected using certain other probability distributions.

Statistically self-similar subsets of any self-similar set may be constructed using a similar percolation process.
Let $\{f_1,\ldots,f_m\}$ be an IFS on $\mathbb{R}^n$ given by \eqref{simform} and let $E_0$ be its attractor. Percolation on $E_0$ may be performed by retaining or deleting components of the natural hierarchical construction of $E$ in a  random but self-similar manner. Let $0<p<1$ and let $D \subset \mathbb{R}^n$ be a non-empty compact set such that $f_i(D) \subset D$ for all $i$. We select a subfamily of the sets $\{f_1(D), \ldots, f_m(D)\}$ where each $f_i(D)$ is selected independently with probability $p$ and write $D_1$ for the union of the selected sets. Then, for each selected $f_i(D)$, we choose sets from
$\{f_{i}f_{1}(D), \ldots, f_{i}f_{m}(D)\}$ independently with probability $p$  independently for each $i$, with the union of these sets comprising $D_2$. Continuing in this way, we get a nested hierarchy $D \supset D_1\supset D_2 \supset \cdots$ of random compact sets, where $D_k$ is the union of the components remaining at the $k$th stage. The random percolation set is $E  = \cap_{k=1}^\infty D_k$, see Figure 3. When the underlying IFS has dense rotations, Falconer and Jin \cite{FJ} extended the ergodic theoretic methods of \cite{HS} to random cascade measures to obtain a random analogue of Theorem \ref{denserots}.

 \begin{figure}[h]
\begin{center}
\includegraphics[scale=0.22]{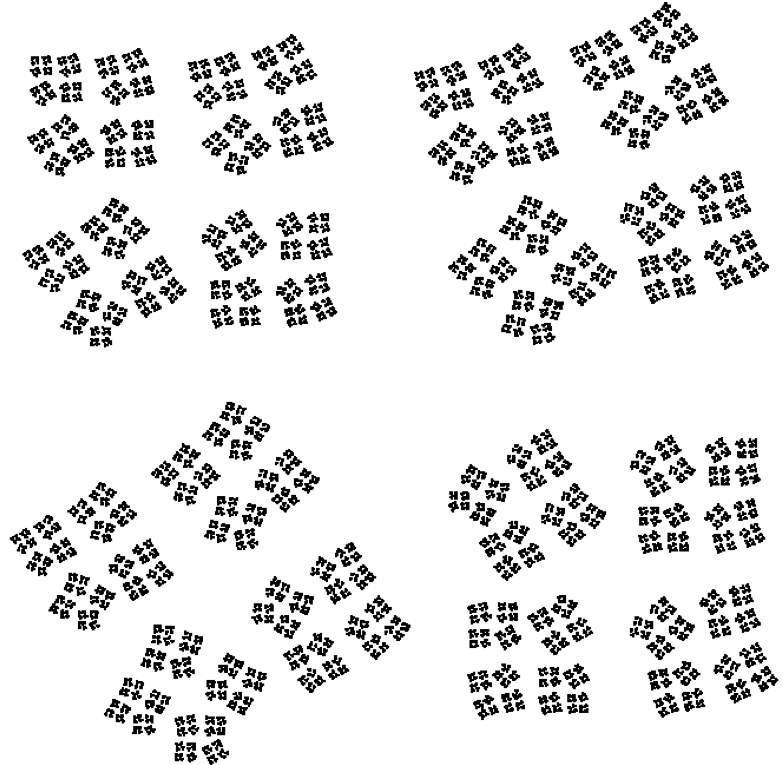}
\qquad\qquad
\includegraphics[scale=0.22]{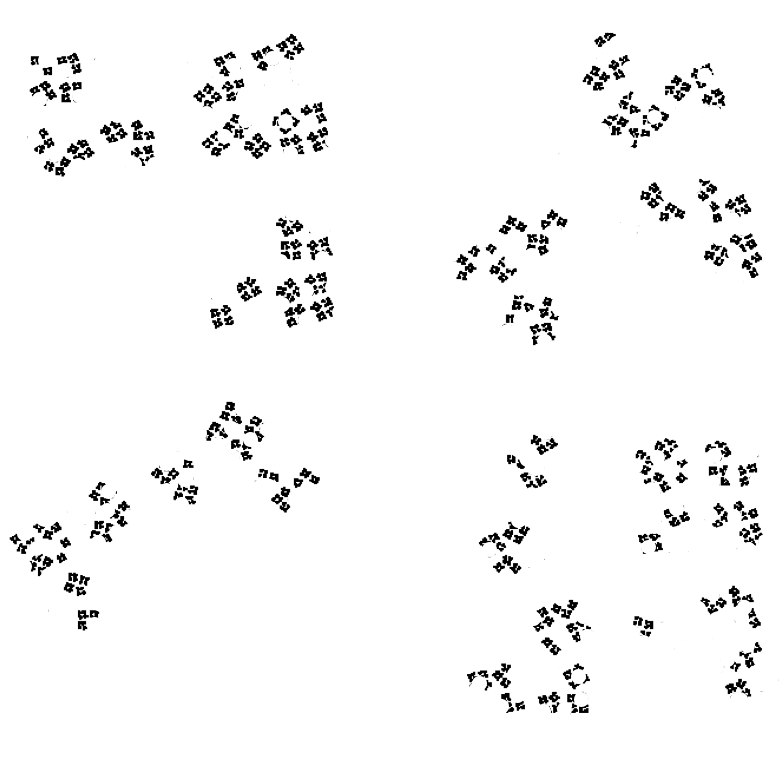}

Figure 3: A self-similar set with dense rotations and a subset obtained by the percolation process.
\end{center}
\end{figure}

\begin{theo}{\rm \cite{FJ}}\label{denserotsperc}
Let $E_0\subset \mathbb{R}^n$ be a self-similar set with dense rotation group and let  $E\subset E_0$ be the percolation set described above. 
If $p>1/m$ there is positive probability that $E\neq \emptyset$, conditional on which,
$$
\dh\mbox{\rm proj}_V E = \min\{\dh E, m\} \mbox{ for  all } V \in G(n,m).
$$
More generally,  conditional on $E\neq \emptyset$, $\dh g(E) = \min\{\dh E, m\}$ for all $C^1$ mappings $g: E\to \mathbb{R}^m$ without singular  points.
\end{theo}

Recently, Shmerkin and Suomala \cite{SSou} have introduced a very general theory showing that for a class of random measures, termed spatially independent martingales, very strong results hold for dimensions of projections and sections of the measures, and thus of underlying sets, with the conclusions holding almost surely  for projections in all directions or onto all subspaces. Such conclusions are obtained by showing that almost surely the total measures of intersections of the random measures with parameterized deterministic families of measures are absolutely continuous with respect to the parameter. 
Spatially independent measures include measures based on fractal percolation, random cascades and random cutout models.

\section{Further variations and applications of projections}
\setcounter{equation}{0}
\setcounter{theo}{0}
\setcounter{figure}{0}

This discussion has covered just a few of the numerous results which may be traced back to Marstrand's pioneering work. We end with an even briefer mention of some further applications, with one or two references indicating where further information may be found.
\medskip

\noindent {\it Visible parts of sets.} The {\it visible part} $\mbox{Vis}_\theta E$ of a compact set $E\subset \mathbb{R}^2$ from direction $\theta$ is
the set of $x \in E$  such that the half-line from $x$ in direction $\theta$ intersects $E$ in the single point
$x$; thus $\mbox{Vis}_\theta E$ may be thought of as the part of $E$ that can be `seen from infinity' in direction $\theta$. It is immediate from Marstrand's Theorem \ref{marthm} that, for almost all $\theta$, 
$$\dh  \mbox{Vis}_\theta E = \dh E \mbox{ if } \dh E \leq1 \quad \mbox{ and } \quad \dh \mbox{Vis}_\theta E \geq 1 \mbox{ if } \dh E \geq 1.$$
 It has been conjectured that if $\dh E \geq 1$  then $\dh  \mbox{Vis}_\theta E = 1$ for almost all $\theta$, but so far this has only been established for rather specific classes of $E$. The conjecture is easily verified if $E$ is the graph of a function (the only exceptional direction being perpendicular to the $x$-axis), see \cite{JJMO}. It is also known for quasi-circles \cite{JJMO} and for Mandelbrot percolation sets \cite{AJJRS}. For self-similar sets, it holds if $E$ is connected and the rotation group is finite  \cite{AJJRS}, and also if $E$ satisfies the open set condition for a convex open set such that $\proj_\theta E$ is an interval for all $\theta$ \cite{FF} (in this case $E$ need not be connected). The analogous conjecture in higher dimensions, that the dimension of the visible part of a set $E\subset \mathbb{R}^n$ equals $\min\{\dh E, n-1\}$, is also unresolved if $\dh E> n-1$.

\medskip

\noindent {\it Projections in infinite dimensional spaces.} Infinite-dimensional dynamical systems may have finite dimensional attractors. When they are studied experimentally what is observed is essentially a projection or `embedding' of the attractor into Euclidean space and infinite-dimensional versions of the projection theorems can relate these projections to the original attractor. Let $E$ be a compact subset of a Banach space $X$ with box-counting (or Minkowski) dimension $d$. Hunt and Kaloshin \cite{HK} show that for almost every projection or bounded linear function $\pi: X \to \mathbb{R}^m$ such that $m >2d$,
$$\frac{m-2d}{m(1+d)} \dh E \leq \dh \pi(E) \leq \dh E.$$
Here `almost every' is interpreted in the sense of {\it prevalence}, which is a measure-theoretic way of defining sparse and full sets for infinite-dimensional spaces. The book by Robinson \cite{Rob} provides a recent treatment of this important area.
\medskip

\noindent {\it Projections in Heisenberg groups.}   The Heisenberg group $\mathbb{H}^n$ is the connected and simply connected nilpotent Lie group of step 2 and dimension $2n+1$ with 1-dimensional
center, which may be identified topologically with $\mathbb{R}^{2n+1}$.  However, the Heisenberg metric $d_H$, which is invariant under the group action, is very different from the Euclidean metric, and in particular the Hausdorff dimension of subsets of $\mathbb{H}^n$ depends on which metric is used in the definition. Despite the lack of isotropy, there is enough geometric structure to enable families of projections to be defined, and it is possible to get bounds for the dimensions of certain projections of a Borel set $E$ in terms of the dimension of $E$, where the dimensions are defined with respect to $d_H$, see \cite{BDFMT, BFMT,Mat5}. 

\medskip

\noindent {\it Sections of sets.} Dimensions of sections or slices of sets, which go hand in hand with dimensions of projections, also featured in Marstrand's 1954 paper \cite{Mar}. He showed essentially that, if $E \subset \mathbb{R}^2$ is a Borel or analytic set  of Hausdorff dimension $s>1$, then for almost all directions $\theta$, $\dh {\proj_\theta^{-1} x} \leq s-1$ for almost  all $x \in V_\theta$, with equality for a set of $x \in V_\theta$ of positive Lebesgue measure. Here $\proj_\theta: \mathbb{R}^2 \to   V_\theta$ is orthogonal projection onto $V_\theta$, the line in direction $\theta$. The natural higher dimensional analogues were obtained by Mattila \cite{Mat, Mat1, Mat4} using potential theoretic arguments. Most of the aspects discussed above for projections have been investigated for sections, including packing dimensions \cite{FM,JM}, exceptional directions \cite{Orp4}, self-similar sets \cite{FJ1,Fur} and fractal percolation sets \cite{FJ1,SSou}.

\medskip

\noindent {\it Projections of measures.}  
For $\mu$ a Borel measure on $\mathbb{R}^n$ with compact support such that $0< \mu(\mathbb{R}^n) <\infty$, the projection $\proj_V\mu$ of $\mu$ onto a subspace $V \in G(n,m)$ is defined in the natural way, that is by 
$(\proj_V\mu)(A) = \mu\{x \in \mathbb{R}^n: \pi_V(x) \in A\}$ for Borel sets $A$ or equivalently by 
$\int f(t) d(\proj_V\mu)(t) = \int \proj_V(x) d\mu(x)$ for continuous $f$. The support of $\proj_V\mu$ is the projection  onto $V$ of the support of $\mu$, so it is not surprising that many of the results for projection of sets have analogues for projection of measures. Indeed many projection results for sets are obtained by putting a suitable measure on the set and examining projections of the measure, as in Kaufman's proof in Section 2. There are many ways of quantifying the fine structure of measures, and the way these behave under projections have been investigated in many cases.

For example, the {\it lower pointwise} or {\it local dimension} of a Borel probability measure $\mu$ on $\mathbb{R}^n$ at $x\in \mathbb{R}^n$ is given by 
$\underline{\dim}_\mu(x) = \varliminf_{r \to 0} \log\mu(B(x,r))/\log r$, with a corresponding definition taking the upper limit for the  {\it upper pointwise dimension} . Then, for almost all every subspace $V \in G(n,m)$ and $\mu$-almost all $x\in \mathbb{R}^n$,
$$  \underline{\dim}_\mu(\proj_V x) = \min\{ \underline{\dim}_\mu( x),m\} \quad \mbox{and} \quad
 \overline{\dim}_\mu(\proj_V x) = \min\{ \overline{\dim}_\mu( x),m\},$$
see   \cite{FOn,HS,HT, HK1,Zah}.
The {\it $($lower$)$ Hausdorff dimension of a measure} $\mu$ is defined as $\dh \mu =\inf\{\dh A : \mu(A)>0\}$. It follows easily from the projection properties of pointwise dimension that
$$ \dh(\proj_V \mu) = \min\{ \dh \mu( x),m\}.$$
The $L_q$-dimensions of projections are examined in \cite{HK1},  for the multifractal spectrum see \cite{BB,Ols,One}, and for packing dimension aspects see \cite{FM}.

For a special case of projection of measures, let $M$ be a compact Riemann surface and $\proj:T^1M \to M$ be the natural projection from the unit tangent bundle $T^1M$ to $M$.  Let $\mu$ be a probability measure on $T^1M$ that is invariant under the geodesic flow on  $T^1M$. Ledrappier and Lindenstrauss \cite{LL} showed that if 
$\dh \mu \leq 2$ then $\dh \proj \mu = \dh \mu$, and if $\dh \mu> 2$ then $\proj \mu$ is absolutely continuous. However, the analogous conclusion fails if the base manifold has dimension 3 or more, see \cite{JJL}. 

\section{Conclusion}
\setcounter{equation}{0}
\setcounter{theo}{0}
\setcounter{figure}{0}

If this article does nothing else, it should demonstrate just how much of fractal geometry has its roots in Marstrand's 1954 paper. If further evidence is needed, there are hundreds of citations of the paper in Math Sci Net and Google Scholar, despite these indexes only including relatively recent references.

This survey of projection results has been brief and far from exhaustive and there are many more related papers. For a both broader and more detailed coverage of various aspects of projections, the books by Falconer \cite{Fal3,Fal} and Mattila \cite{Mat1,Mat2} and the survey articles by Mattila \cite{Mat3,Mat4,Mat5} may be helpful.

\end{document}